\documentclass{elsarticle}
\usepackage[utf8]{inputenc}

\usepackage{amssymb}
\usepackage{amsmath}
\usepackage{subcaption}

\usepackage{natbib}
\usepackage{graphicx}
\usepackage{hyperref}
\usepackage{bm}

\usepackage{algorithm}
\usepackage{algpseudocode}

\journal{Computers \& Fluids}

\begin{document}

\begin{frontmatter}

\title{Accelerating multigrid solver with generative super-resolution}

\author[1,2]{Francisco Holguin}
\ead{fholgui1@jhu.edu}

\author[3]{GS Sidharth}

\author[4]{Gavin Portwood}

\affiliation[1]{organization={Johns Hopkins Applied Physics Lab},
                city={Laurel},
                postcode={20723},
                state={MD},
                country={USA}}
% addressline={}

\affiliation[2]{organization={Whitling School of Engineering, Johns Hopkins University},
                city={Baltimore},
                postcode={21218},
                state={MD},
                country={USA}}

%\affiliation[3]{organization={Los Alamos National Laboratory},
%                city={Los Alamos},
%                state={NM},
%                postcode={87545},
%                country={USA}}
                
\affiliation[3]{organization={Dept. of Aerospace Engineering, Iowa State University},
                city={Ames},
                state={IA},
                postcode={50011},
                country={USA}}

\affiliation[4]{organization={Lawrence Livermore National Laboratory},
                city={Livermore},
                state={CA},
                postcode={94550},
                country={USA}}

\begin{abstract}
The geometric multigrid algorithm is an efficient numerical method for solving a variety of elliptic partial differential equations (PDEs). The method damps errors at progressively finer grid scales, resulting in faster convergence compared to iterative methods such as Gauss-Seidel. The prolongation or coarse-to-fine interpolation operator within the multigrid algorithm, lends itself to a data-driven treatment with deep learning super-resolution, commonly used to increase the resolution of images. We (i) propose the integration of a super-resolution generative adversarial network (GAN) model with the multigrid algorithm as the prolongation operator and (ii) show that the GAN-interpolation can improve the convergence properties of multigrid in comparison to cubic spline interpolation on a class of multiscale PDEs typically solved in fluid mechanics and engineering simulations. We also highlight the importance of characterizing hybrid (machine learning/traditional) algorithm parameters.
\end{abstract}

\begin{keyword}
%% keywords here, in the form: keyword \sep keyword

%% PACS codes here, in the form: \PACS code \sep code

%% MSC codes here, in the form: \MSC code \sep code
%% or \MSC[2008] code \sep code (2000 is the default)

\end{keyword}

\end{frontmatter}

\section{Introduction}
Partial differential equations describe the behavior of many physical phenomena, from the smallest quantum scales to the largest structures in the universe. Efficient numerical methods are crucial to simulate these complex phenomena with current and foreseeable computational resources. The application of machine learning (ML) tools to accelerate and improve the accuracy of physical simulations is an increasingly active area of research, leveraging the surge of advancements in machine learning architectures over the last decade \cite[e.g.][]{wei2018, peurifoy2018, raissi2019physics,ranade2020,sanchez2020}.

While replacing conventional numerical solvers with machine learning models is certainly promising, improving heuristic operators within existing  mathematically-formally derived numerical methods can allow for interpretable computational gains, easier implementation and more rapid yet controlled deployment\citep[][]{markidis2021old}.

We focus on geometric multigrid methods, which are effective at solving multiscale PDEs used in physics and engineering \citep[e.g.][]{trottenberg2001, gs2020multiscale}. The multiscale nature of multigrid methods motivates the integration of multiscale ML models, in particular deep learning methods \citep{lecunhinton2015}. Recent investigations \citep[e.g.][]{markidis2021old,margenberg2022} similar to ours have explored the use of deep neural networks in a geometric multigrid solver by evaluating the effect of network properties on solution accuracy, runtime, and training time. They find that deep learning applications to multigrid can plausibly accelerate the algorithms while still accurately recovering the true solution. Super-resolution (SR) deep learning models, in particular, are multiscale in nature and have the popular use of transforming photographs from low to high resolution while retaining a visually realistic result.

The goal of our work is to explore the integration of a deep learning super-resolution model into the geometric multigrid algorithm. We focus on applying the algorithm to the two-dimensional fluid system described by the pressure-Poisson formulation of the Navier-Stokes equations. We discuss the conversion of a super-resolution model from a photographic to a fluid system and the integration of the ML model within the algorithm. We also investigate the algorithm performance, including scaling and convergence, for several key parameters describing the combined ML-multigrid algorithm.

\section{Methods}

\subsection{The Pressure-Poisson equation in the incompressible Navier-Stokes equations}
%brief description of the system to give context
The incompressible Navier-Stokes (NS) equations govern fluid flow. The equations have a quadratic non-linearity that leads to multiscale solutions (such as in the case of hydrodynamic instabilities and turbulence) and represent fluid physics in widespread applications in engineering and scientific phenomena. The incompressible NS equations \citep{quartapelle2013numerical} are given by

\begin{equation}
\begin{split}
\frac{\partial \bf{u}}{\partial \rm{t}} + (\bf{u} \cdot \nabla) \bf{u} &= - \nabla p + \nabla \sigma \\
\nabla \cdot \bf{u} &= 0
\end{split}
\label{ns_eq}
\end{equation}

where $u$ is the fluid velocity and $p$ is the pressure.
%We have ignored viscosity and body forces. 
This equation can be discretized and evolved in time for the fluid velocity and pressure.
The role of the pressure in the incompressible equations is to maintain mass conservation or zero velocity divergence.
It is possible to take the divergence of the momentum equation and rewrite Eq. \ref{ns_eq} in the form

\begin{equation}
\begin{split}
&\frac{\partial \bf{u}}{\partial \rm{t}} + (\bf{u} \cdot \nabla) \bf{u} =  -\nabla p \\
& \nabla^2 p = - \nabla \cdot (\bf{u} \cdot \nabla)\bf{u},
\end{split}
\end{equation}
the discretization of which ensures zero production of velocity divergence in the solution.
As a consequence, the fluid pressure $p$ is now described by a Poisson equation. Time-evolution of this system generally involves a splitting that solving the Poisson equation for the pressure given a velocity field increment due to convection-diffusion terms, and then using the pressure correction to evolve the velocity field. The details of the boundary conditions are specific to each problem and are not explicitly described here. The Poisson equation for pressure is a non-local operation and is often the rate-determining step in the solution procedure for incompressible convection-dominated flow problems. Accelerating the process of solving the Poisson equation, therefore, directly impacts the speedup in the overall time-integration algorithm. 

%It must be noted that our work restricts to 

%\begin{equation}
%    \nabla^2 p =  f(x,y),
%\label{poissoneq}
%\end{equation}
%where the solution $p$ is discretized on a grid, and $f$ is some spatially dependent function, which in this case is a function of the fluid velocity. 

\subsection{The Multigrid method for the Poisson equation}
\label{multigrid_section}

Multigrid methods are particularly effective for solving multiscale differential equations relevant to physics and engineering \citep[e.g.][]{trottenberg2001,gs2020multiscale}. Multigrid methods formally scale linearly with the number of unknowns in a discretized physical system, which makes them attractive to use, as high-performance computing increasingly enables higher-fidelity and larger-scale computations. Generally, the favorable scaling of multigrid methods is due to their ability to iteratively reduce solution errors at multiple scales by obtaining coarse solutions, then utilizing the coarsened solutions to solve increasingly higher-fidelity solutions.  The utilization of coarse fidelity solutions to accelerate finer fidelity solutions is enabled by a prolongation (or interpolation) operator.

Prolongation within a multigrid algorithm lends itself to a data-driven treatment due to the ability of ML techniques to perform data-informed interpolation, as seen frequently seen in super-resolution methods as used in image analysis \citep[c.f.][]{farsiu2004advances,park2003super}. In this work, we consider the applications of enhancing multigrid methods by using a data-driven prolongation operator akin to super-resolution operators used in image analysis. We consider a common Poisson differential equation, frequently written as
\begin{equation}
    \nabla^2 p = \frac{\partial^2 p}{\partial^2 x} + \frac{\partial^2 p}{\partial^2 y} = f(x,y),
\label{poissoneq}
\end{equation}
where the solution $p$ is discretized on a grid, and $f$ is the spatially dependent source term, which is a non-linear function of the fluid velocity for the pressure-Poisson equation. 

\begin{figure}
    \centering
    \begin{subfigure}[t]{0.8\textwidth}
        \includegraphics[width=\textwidth]
        {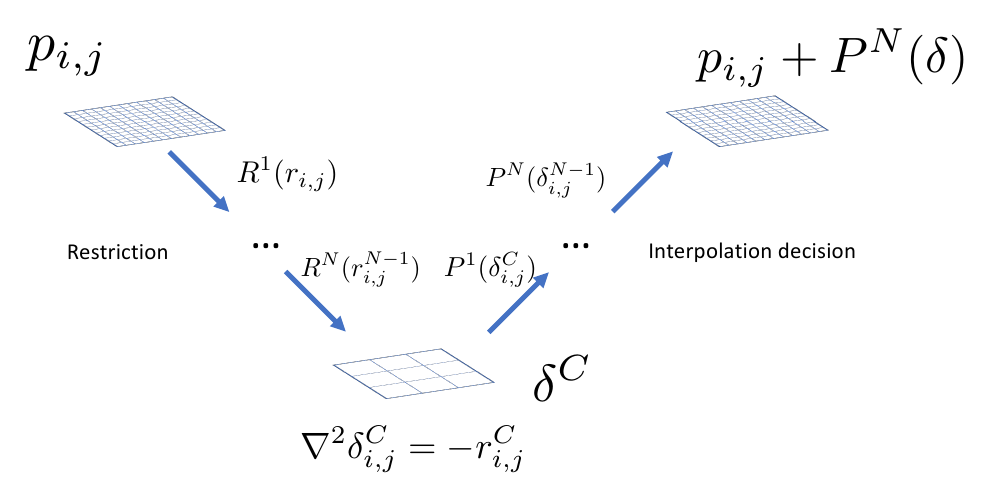}
        \caption{Diagram of a single iteration of a geometric V-cycle multigrid algorithm. The coarsest level is denoted by a superscript `$C$'. Relaxation operations occurs at every resolution level to determine the correction term $\delta^{i}$, although the error correction PDE is only shown for the coarsest level. The restriction and prolongation operators are shown. In our work, there is a decision as to which prolongation operator is used in that multigrid iteration.}
        %{ml_mg_figures_Vcycle_updated2.pdf}
        %\caption{Schematic of a two-level multigrid cycle. Adapted from \citet{chen2001}.}
        \label{vcycle}
    \end{subfigure}
    ~
    \begin{subfigure}[t]{0.8\textwidth}
        \includegraphics[width=\textwidth]{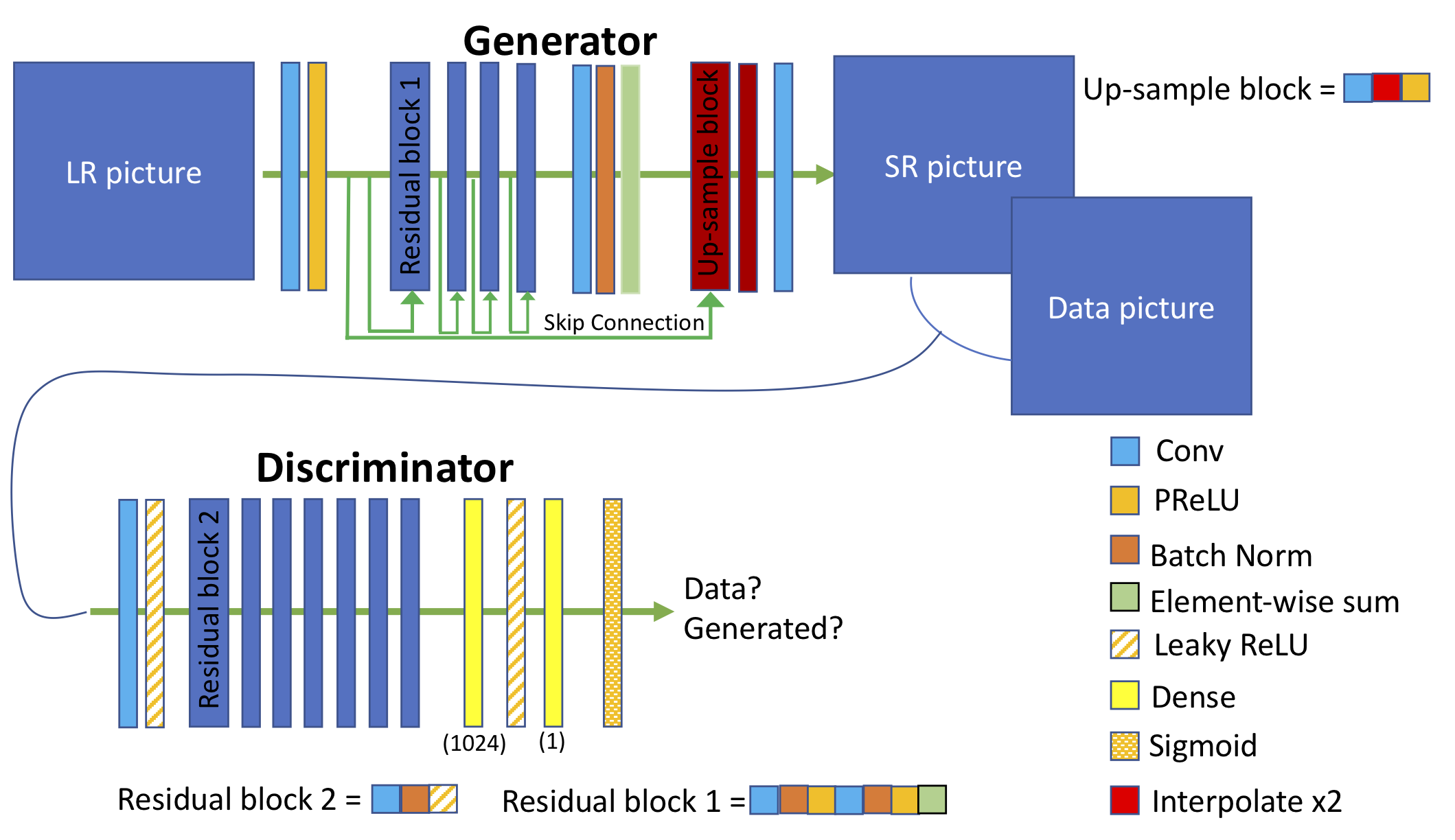}
        \caption{Structure of discriminator and generator we implement, modified from \citet{ledig2017}. }
        \label{gan_structure}
    \end{subfigure}
    \caption{Diagrams for two-level multigrid and GAN interpolation.}
\end{figure}

Figure \ref{vcycle} shows a schematic of the two-level multigrid, where the top horizontal level represents the highest resolution grid $h_2$. Applying the Poisson operator $\nabla^2 p - f$ to the grid $h_2$, using a relaxation method such as Gauss-Seidel, results in an residual $r_f$ relative to the true solution. For these relaxation methods, high-frequency errors (relative to the grid scale) are quickly reduced, while low spatial frequency errors remain. It can be shown \citep[][]{trottenberg2001} that the residual also satisfies the Poisson equation, so we recursively solve the Poisson equation for this error. In order to efficiently suppress low-frequency errors, we apply a restriction operator $\bf{R}$ to coarsen the grid so that previously low spatial frequency errors on the fine grid become high frequency errors on the coarse grid $h_1$. 

Some applications of deep learning networks directly employ surrogates for the PDE, either on the course grid \citep[i.e.][]{markidis2021old} or at each multigrid level \cite{zhang2022hybrid}. We only employ Gauss-Seidel relaxation for these operations. Once the coarsest level is reached, we apply a prolongation operator $\bf{P}$ to the error term in order to correct the solution on the next finer grid. In a uniformly spaced grid, the restriction operator $\bf{R}$ can be as simple as taking only the black points in red-black ordering. The prolongation operator $\bf{P}$ requires more work, as it involves an interpolation scheme to create a higher resolution grid. For example, this interpolation scheme can be a two-dimensional mesh interpolation, such as the bivariate spline function in the SciPy Python library \citep[][]{virtanen2020}. 

\subsection{The GAN model for super-resolution}
\label{GAN}
Deep learning \citep{lecunhinton2015} methods have been able to capture and reproduce the structure of extremely complex physical systems \citep[e.g.][]{portwood21,wang20,he2019}. Even without prior knowledge of the system, the network can learn the behavior and patterns relatively quickly with adequate data.

Traditional supervised learning, such as for classification tasks, involves a network learning a decision boundary via labeled data (e.g. images \citep[][]{al2017review} or text \citep[][]{minaee2021deep}). On the other hand, unsupervised generative models learn a probability distribution from the data, allowing us to then sample from the distribution, and create new synthetic data. More precisely, as described in \citet{harshvardhan2020comprehensive}, given input $X$ and output $Y$ variables, classification models learn the probability $P(Y|X)$, the probability of $Y$ given a particular input example $X$. Generative models instead learn $P(X|Y)$, the probability distribution of inputs (or features) corresponding to the output, and $P(Y)$, the probability of generating a particular synthetic output. By learning the probability distribution of features that correspond to outputs, the model can generate realistic synthetic data.

Super-resolution \citep[e.g.][]{ledig2017, yang2019deep} in deep learning refers to the application of generative networks to low-resolution images in order to predict a high-resolution version. There are many generative model architectures, such as hidden Markov models, Boltzmann machines, variational autoencoders, normalizing flows, and generative adversarial networks \citep[][]{harshvardhan2020comprehensive, kobyzev2020normalizing}. A Generative Adversarial Network (GAN) \citep{goodfellow2014} is particularly effective as an accurate, super-resolution model. GANs do not explicitly learn the underlying data probability distribution: through their adversarial structure, they develop an approximate, efficient lower dimensional latent representation \citep[][]{creswell2018generative,jabbar2021survey}.

The GAN network consists of a generative (`Generator' $\bf{G}$) model with a discriminative (`Discriminator' $\bf{D}$) model. The idea is to train the Generator to produce increasingly more realistic data, while the Discriminator trains to distinguish between real and synthetic samples. Training occurs as follows: (i) use a fixed Generator to generate a synthetic data sample, (ii) train the Discriminator via labeled (real or synthetic) data samples with categorical loss $l_{\rm{Dis}}$, and (iii) fix the Discriminator and train the Generator by inputting new synthetic and real data into the Discriminator and constructing an appropriate loss function $l^{SR}$ that optimizes both for the distance of the synthetic data to the real data and against the performance of the discriminator. We implement the losses described in \citep[][]{ledig2017}:
\begin{equation}
\begin{split}
l_{\rm{Dis}} = \ & l_{\rm{BCE}} \\
l^{SR} = \ & l^{SR}_{MSE} + 10^{-3} l^{SR}_{Gen} \\
l^{SR}_{MSE} = \ & \frac{1}{r^2 W H} \sum_{x=1}^{rW} \sum_{y=1}^{rH} (I_{x,y}^{HR} - G_{\theta} (I^{LR})_{x,y} ) ^2 \\
l^{SR}_{Gen}  = \ & \sum_{n=1}^N - \rm{log} D_{\theta} (G_{\theta} (I^{LR})),
\end{split}
\end{equation}
where $I^{LR}$ and $I^{HR}$ are the low/high resolution grids described by a tensor of size $W \times H \times C$ and $rW \times rH \times C$  respectively, $r$ is the up-sampling factor along each spatial dimension, $\theta$ represents the model parameters in a particular network, and $l_{\rm{BCE}}$ is the binary-cross entropy loss. 

In our work, we use the GAN architecture from \citet{ledig2017}, following the example in \citet{birla2018}, as shown in Figure \ref{gan_structure}. Unlike \citet{ledig2017} who use the pre-trained VGG loss designed for photographic image recognition, we use mean squared error (MSE) loss function for the generator. We keep the same dimension of `color' channels $C=3$, but set each channel to have equal data. This choice allows minimal architecture changes from the original image focused-network, as well as flexibility for the inclusion additional information in future work.

\section{Multigrid with GAN super-resolution}

\subsection{Dividing the domain and GAN interpolation}
\begin{figure}
    \centering
    \includegraphics[width=0.9\textwidth]{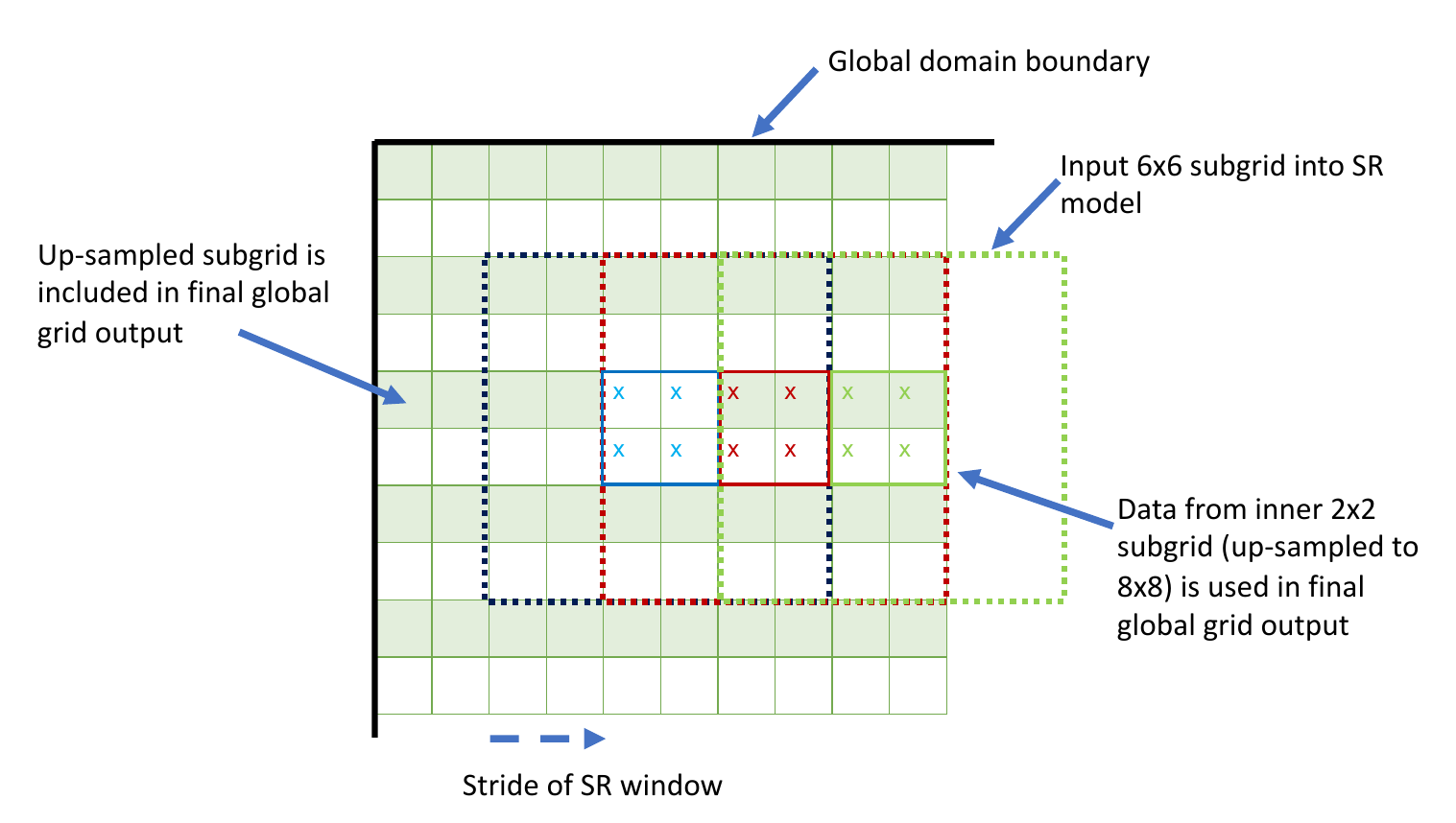}
    \caption{Schematic of ghost cell-like splitting of the domain. The blue dotted line shows the input grid of size $n_s \times n_s$ for the SR GAN. Only the upscaled data from the inner $2\times 2$ grid is kept for the final grid reconstruction.  The next input grid is denoted by red dotted line, which is shifted 2 grid elements to the right, so that the inner $2\times 2$ grid data (marked by `x') that is kept is next to the previous upscaled data. The processes continues with the next grid in green, and so on until. The algorithm proceeds similarly in downwards. }
    \label{domain_split_ex}
\end{figure}

The GAN maps a ($n_s \times n_s$) (here $n_s=6$) subset of the discretized domain to a ($n_l \times n_l$) grid. These windowed inputs allow for spatially-invariant training and prediction. A small value of $n_s$ also allows for generation of more training samples from limited datasets. Prediction for the whole domain involves stitching together the windowed outputs. Ghost cells are implemented by splitting the input domain into overlapping windows and utilizing only subset of that particular window's error prediction for the global high resolution grid. Figure \ref{domain_split_ex} shows a schematic for how this domain splitting and data reconstruction occurs. Only the upscaled data from the grid marked by `x' is used in the reconstruction of the final high resolution grid. This implementation reduces numerical artifacts in the SR high resolution grid as the adjacent individual tiles are predicted by the SR GAN operator with overlapping data as input. Before input in the GAN, we normalize the grid data in log space to the range (-1,1). Unlike photographic image cells, which typically span the range (0,255), pressure data spans many orders of magnitude. Operating on the log of the data greatly improves training of the GAN.

The GAN prolongation operator works as follows (see \ref{alg_appendix} for the algorithm description):

\begin{itemize}
    \item Normalize the grid values from $p_{i,j}$ to $q_{i,j}$, by using a symmetric log function $T : \mathbb{R} \to \Omega \in [-1,1]$. This function is constructed using another function $T': \mathbb{R} \to \Omega \in [0,1]$, so that $T(p_{i,j}) = T'(|p_{i,j}|) \ \rm{sgn}(p_{i,j})$. The function $T'$ linearly maps the log of the absolute values. We define the min and max magnitudes ${p}_{\textrm{min}}$ and ${p}_{\textrm{max}}$ to encompass the values expected in the data. The resulting normalization is such that $T({p}_{\textrm{min}} ) ) = 0 $ and $T({p}_{\textrm{max}} ) ) = 1 $, and $T({p}_{\textrm{min}} ) ) = 0 $ and $T({p}_{\textrm{max}} ) ) = -1 $  .
    \item Divide the normalized coarse grid $q_{i,j}$ into overlapping $n_s \times n_s$ windows (i.e. assuming a $n_s^2$ kernel with stride = 2). Apply GAN interpolation to produce a set of high resolution $n_l^2$ window kernels. Here $n_l=24$.
    \item Construct the fine grid by assigning the central $(n_s+2)^2$ fine window values to the final fine grid. Any non-overlapping areas (i.e the edges at the global boundary) are assigned to the output window  in that region.
    \item Transform and scale back to original, un-normalized pressure values.
\end{itemize}

\subsection{Data and training}
\label{training}

The pressure-Poisson formulation of the incompressible NS equations is a physical system with the form of Eq. \ref{poissoneq}. The pressure-Poisson equation is used to ensure mass conservation in incompressible flows and is derived from taking the divergence of the fluid momentum equation 
\begin{equation}
    \nabla^2 p = \nabla \cdot (\nu \nabla^2 u  - ( u \cdot \nabla ) u    )  = f(x,y)
    \label{poisson_ns}
\end{equation}

Here $\nu$ is the fluid viscosity, $u$ is the fluid velocity, and $p$ is the fluid pressure. Simulating the system in time involves calculating the source term $f(x, y)$ from the velocity field $u$, solving the Poisson equation for pressure (e.g. using a multigrid method), and then using that pressure to update the velocity field. We produce a truth dataset of 200 pressure fields on a grid size $192^2$ by evolving a biperiodic velocity field from a broadband initial condition and Reynolds number $\rm{Re}_{\rm{train}} = 2000$ using the incompressible NS equation. 

We construct the training set of 1000 pressure grids as follows:
\begin{itemize}
    \item Select one of the 200 pressure fields randomly and restrict (downsample) the grid by a random power $l$ of 2 on each grid length, such that the smallest grid size is $(n_s=6)^2$. The result is a grid of size $(192 / 2^l)^2 $.
    \item From the grid above, pick a random $6^2$ window. The total size of the space of such windows is $17,000$.
    \item Transform and rescale the window data range from $\displaystyle [- \infty, \infty]$ to $\displaystyle [- 1, 1]$ in the same manner as described in the GAN prolongation method in Section \ref{GAN}, using the global min and max of the truth dataset. We used $  |p_{\rm{min}}| = 10^{-10}$ and $|p_{\rm{max}}| = 10^{-3}$.
\end{itemize}

\begin{figure}
    \centering
    \includegraphics[width=\textwidth]{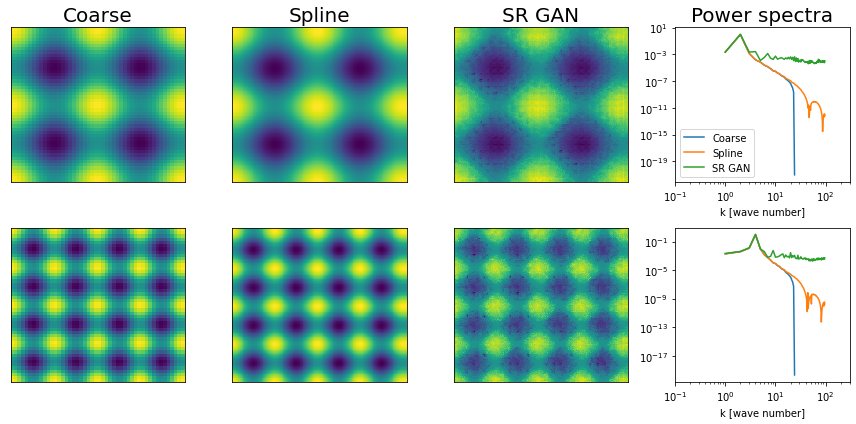}
    \caption{Test of super-resolution GAN model on pressure data generated with a sinusoidal function at a fixed length scale. The sinusoidal function is $p(x,y) = \rm{cos}(nx) + \rm{cos}(ny)$, where ($x,y$) are the spatial coordinates in the range $[0,2\pi]$ and $n$ is an integer. The first and second rows correspond to $n=2$ and $n=4$ respectively. The first column shows generated data at $N=48$. The second and third columns show the interpolated high-resolution grid at $N=192$ using a spline and SR GAN method, respectively. The fourth column shows the power spectra of spatial scales for all three grids, each normalized to their peak value. Please note that the slope of the SR GAN spectra is not zero, and appears small only due to the y-axis range. }
    \label{SR_GAN_test}
\end{figure}

In Figure \ref{SR_GAN_test}, we examine the behavior of the trained SR GAN model by using the model to interpolate simple, single-mode test data. We generate these test pressure grids with the function  $p(x,y) = \rm{cos}(nx) + \rm{cos}(ny)$, where $(x,y) \in [0,2\pi]$ are the spatial coordinates and $n$ is an integer. We show examples for large and small length scale data (rows 1 and 2). Broadly, we see that the SR GAN operator produces a reasonable interpolation. A comparison of the exact norm difference between the spline and SR GAN grids is not a useful comparison because the generated data happens to contain no smaller length scales (i.e. we generated it with a sinusoidal function that is locally linear). The purpose of this test is to demonstrate how the high-resolution predictions differ in terms of spatial scales. The last column shows the power spectra of the coarse, spline interpolated, and SR GAN interpolated grids. Both the spline and SR GAN methods capture the large spatial scales (low $k$) similarly, as shown by their overlapping power spectra with that of the coarse grid. The spline interpolated grid spectra follow the spectrum of the coarse grid and continue the same power law below the coarse scale. In contrast, the SR GAN interpolated grid power spectra is higher at smaller scales as the power law is shallower. This behavior is expected, as the spline operator is simply extending the coarse data based on local data points. The SR GAN operator predicts what the finer structure should be based on the learned representation of the fluid pressure data.

\begin{figure}
    \centering
    \includegraphics[width=0.75\textwidth]{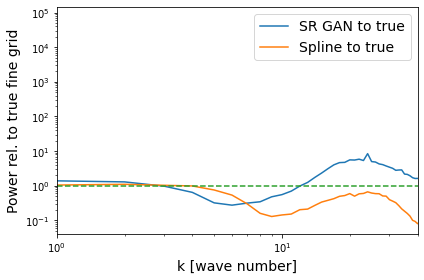}
    \caption{Examination of super-resolution GAN model on grids from the test set. The predicted power spectra are smoothed with $\rm{N}_{\rm{smooth}}=5$ via the Gauss-Seidel iterations. Each line is the power spectra relative to the true fine solution, averaged from 50 test grids. The dashed line shows the line value of one.}
    \label{SR_GAN_test_set}
\end{figure}

Figure \ref{SR_GAN_test_set} is similar to Figure \ref{SR_GAN_test} except that we applied the SR GAN and spline interpolation operators to test data . We smooth the predicted grids in order to better match properties during the multigrid algorithm, as there is inter-iteration smoothing. Both interpolation operators predict large spatial scales similarly well. Spline interpolation tends to under-predict smaller spatial scales, while the SR GAN interpolation tends to overpredict these scales. The operators are the closest to the true grid power spectra at different scales. These differences in grid predictions manifest as differences in the number of iterations to convergence as we discuss in Section \ref{results}. 

\subsection{Algorithm application}
\begin{table}[]
\begin{tabular}{llll}
 & Parameter            & Fiducial value                            & Range of values tested                                \\ \cline{2-4} 
 & \multicolumn{3}{c}{Multigrid}                                                                                            \\ \cline{2-4} 
 & $N_{\rm{iter}}$          & 300                                       & 300                                                   \\
 & $N_{\rm{grid}}$        & 192                                       & 96,192,384,768                                        \\
 & $r_{\rm{min}}$         & 12                                        & 6,12,24,48                                            \\
 &                      &                                           &                                                       \\
 & $N_{\rm{smooth, pre}}$ & 10                                        & 10                                                    \\
 & $N_{\rm{smooth}}$       & 20                                        & 15, 20, 25                                            \\
 & $N_{\rm{step}}$        & 4                                         & 4                                                     \\
 &                      &                                           &                                                       \\
 &                      &                                           &                                                       \\ \cline{2-4} 
 & \multicolumn{3}{c}{Data}                                                                                                 \\ \cline{2-4} 
 & Training set size    & 1000                                      & 1000                                                  \\
 & Testing set size     & 100                                       & 100                                                   \\
 & $Re_{\rm{train}}$      & 2000                                         & 2000                                                    \\
 &                      &                                           &                                                       \\ \cline{2-4} 
 & \multicolumn{3}{c}{GAN}                                                                                                  \\ \cline{2-4} 
 & Layer parameters     & see \citet{ledig2017}                               &                                                       \\
 & $n_{\rm{up}}$          & 2 (= $N_{\rm{step}}/2)$                       &                                                       \\
 &                      &                                           &                                                       \\ \cline{2-4} 
 & \multicolumn{3}{c}{GAN+Multigrid}                                                                                        \\ \cline{2-4} 
 & $n_{\rm{s}}$           & 6                                         &                                                       \\
 & $n_{\rm{l}}$         & $n_{\rm{s}} n_{\rm{up}}$                       &                                                       \\
 &                      &                                           &                                                       \\
 & $N_{\rm{GAN}}$         & $1/N_{\rm{iter}}$ \ $\rm{or}$ \ $N_{\rm{iter}}$ & $1/N_{\rm{iter}}$,  1/10, 1/5, 1, 5,  10, $N_{\rm{iter}}$ \\
 & $S_{\rm{thres}}$       & 0 or 1                                    & 0, $10^{-5}, 10^{-3}$, 1                                \\
 & $Re_{\rm{test}}$       & 2000                                         & 100, 2000, 5000                                              
\end{tabular}

\footnotesize{\textbf{Parameters} From top to bottom the rows contain: (1) total number of multigrid iterations; (2) side length of square grid; (3) coarsest grid side length; (4) number of GS smoothing iterations on initially random pressure grid; (5) number of GS smoothing iterations in between multigrid iterations; (6) resolution change (change in grid side length) between multigrid levels; (7),(8) number of grids in training and test set respectively; (9) Reynolds number of training grids; (10) GAN model paramters and choices; (11) up-sample factor (change in grid side length) that SR GAN model performs; (12) side length of SR interpolation kernel/window; (13) side length of SR interpolated kernel; (14) number of SR GAN operations per spline operation; (15) switch for changing interpolation operator at a particular stage of convergence; (16) Reynolds number of test grids.}
\caption{}
\label{param_table}
\end{table}

We solve Eq. \ref{poisson_ns} using a test set source grid $f(x,y)_{i}$ from a velocity field $u_{i}$. We set the viscosity $\nu = 0$ for simplicity. Our implementation of this hybrid multigrid algorithm is characterized by many parameters related to both the deep learning and non-deep learning parts of the algorithm. \ref{alg_appendix} contains an algorithm description. Table \ref{param_table} summarizes the parameter choices we make, including the fiducial choices described in Section \ref{motivatingresults}. The multigrid parameters $N_{\rm{smooth, pre}}$ and $N_{\rm{smooth}}$ describe the number of Gauss-Seidel smoothing applications on the initial grid and at the beginning of each multigrid iteration respectively. This additional smoothing is often necessary to ensure smoothness at the grid scale in order to avoid aliasing errors. The parameters $N_{\rm{step}}$ and $r_{\rm{min}}$ control the number of multigrid levels. $N_{\rm{step}}$ is the reduction resolution per multigrid level on each spatial dimension in units of a factor of 2. For example, a restriction or interpolation from a $192^2$ to a $96^2$ grid is $N_{\rm{step}} = 1$, while a restriction from the same starting grid to $48^2$ is $N_{\rm{step}} = 2$. The parameter $r_{\rm{min}}$ denotes the coarsest resolution desired. The pair of parameters $N_{\rm{GAN}}$ and $S_{\rm{thres}}$ each allow for a choice of spline and super resolution interpolation in each individual multigrid iteration. $N_{\rm{GAN}}$ describes the number of GAN interpolations per spline interpolation. For example, $N_{\rm{GAN}}$ = 1 results in alternating spline and GAN interpolations. The edge cases of only spline and only super resolution have $N_{\rm{GAN}} = 1/N_{\rm{iter}} $ and $N_{\rm{iter}}$ respectively. The parameter $S_{\rm{thres}}$ controls which interpolation method is used at a particular integration. It represents a threshold of the norm of the difference between grids $|p_i - p_{i+1}|$ at which to change the interpolation type. This parameter allows us to explore the use of different interpolation methods for different parts of the calculation, such as early-on when the grid is still far from the solution or when the grid is already close to converged.

%
% insert SR MG algorithm here?
% link to apendix?
%

To study the efficacy of a GAN prolongation operator under controlled settings, we implement a two-level V-cycle multigrid algorithm (i.e. fine scale $192^2$ grid directly to a coarse scale $12^2$ grid, corresponding to $N_{\rm{step}}$ = 4 and $r_{\rm{min}} = 12$). A two-level multigrid method with a coarsening factor of 16 is expected to perform poorly, due to the extreme prolongation and interpolation, and therefore serves as an ideal test bed for evaluating the gains made with a super-resolution GAN. We construct $100$ pressure and source term grid pairs for our test set, following the method described in Section \ref{training}, which independent of the training set data. The Reynolds number $\rm{Re}_{\rm{test}}$ for the fiducial test case is the same as the one for the training set, although we also briefly explore grids with different values of $\rm{Re}_{\rm{test}}$, to measure how well the super resolution model performs on data with different scales than it was trained on.

\section{Results}
\label{results}

\subsection{Motivating results}
\label{motivatingresults}

\begin{figure}
    \centering
    \begin{subfigure}[t]{\textwidth}
        \includegraphics[width=\textwidth]{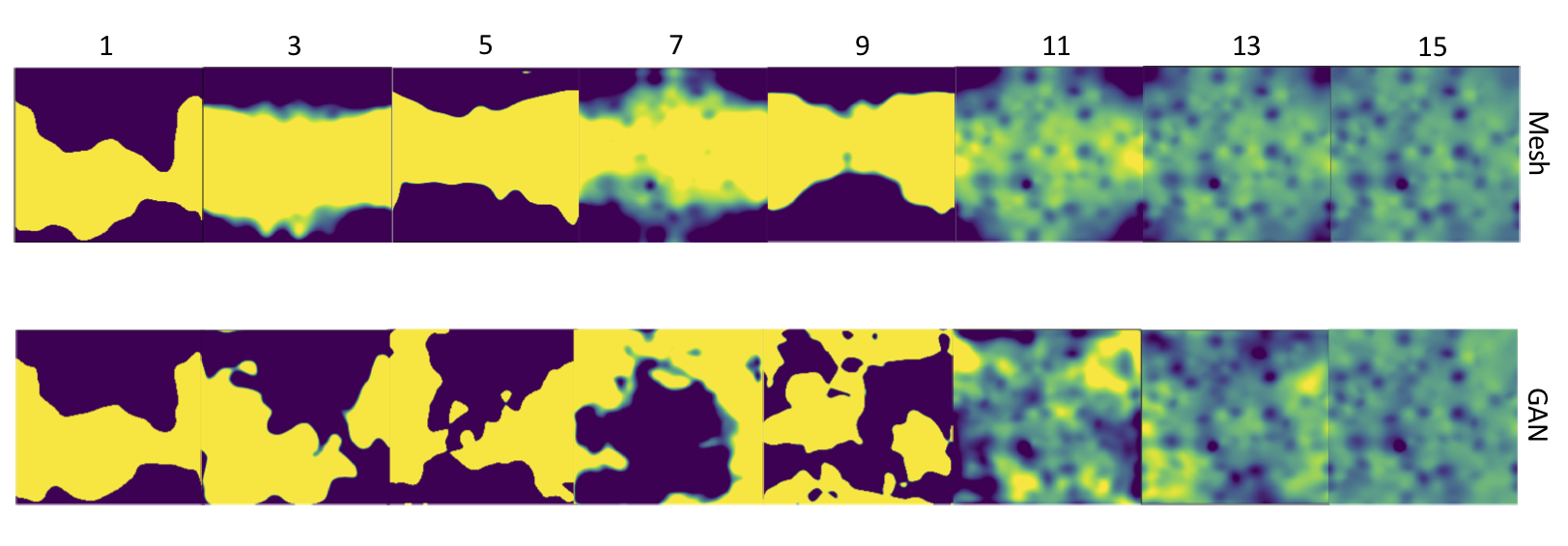}
        \caption{Visual comparison of the grid pressure in the first 15 iterations of solving Eq. \ref{poissoneq} with the multigrid solver starting from an initially random, smoothed grid. The colorscale is fixed to represent the true pressure field. Top: interpolation with traditional mesh interpolation. Bottom: interpolation with GAN super resolution. The GAN version shows higher frequency structure earlier than the spline version.}
        \label{ml_mg_data_series}
    \end{subfigure}
    ~
    \begin{subfigure}[t]{\textwidth}
        \includegraphics[width=\textwidth]{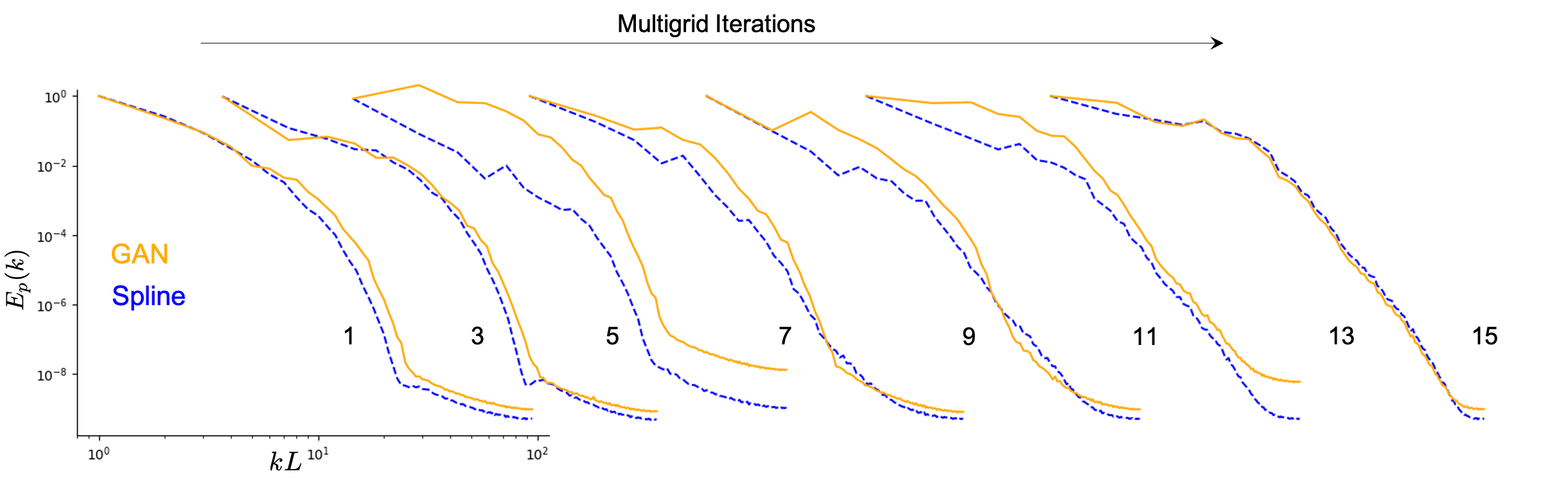}
        \caption{Convergence in $p$ spectral density with increasing multigrid iterations from the same data as Figure \ref{ml_mg_data_series}. As in the visual comparison, the GAN multigrid has higher frequency structure at earlier times. Both multigrid types converge to the same solution by iteration 15.}
        \label{iter}
    \end{subfigure}

\end{figure}

As an initial proof-of-concept test, we solve a two-level multigrid system for 100 grids from the test set and assume the fiducial parameters. We compare two implementations of the interpolation step of the multigrid algorithm, the bivariate spline and the super resolution-based operators. This subsection is based on our previous work in \citet[][]{holguin2021multigrid}.

%%%%

\begin{figure}
    \centering
    \includegraphics[width=\textwidth]{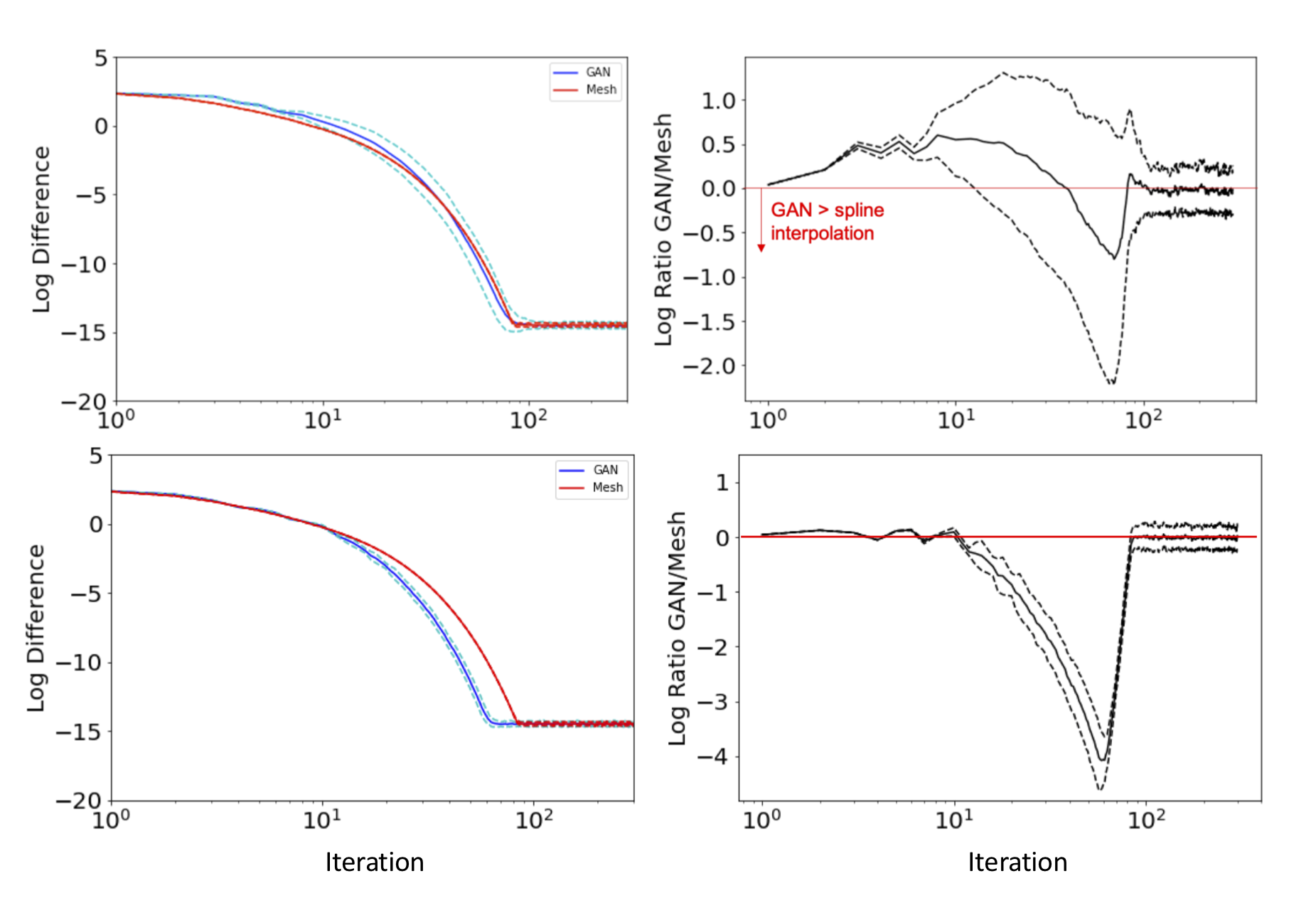}
    \caption{Averaged results of multigrid algorithm solving $100$ different grids. Norm of the difference of grid between iterations as a function of iteration, shown for two choices of prolongation operator, two-dimensional spline and GAN super resolution interpolation. The multigrid parameters used are $N_{\rm{smooth, pre}}$ = 10, $N_{\rm{smooth}} $ = 20,  $N_{\rm{step}}$ = 4, and $r_{\rm{min}}$ = 12. Top: The GAN operator is used for interpolation in every iteration. Bottom: The GAN and spline operators are alternated every other iteration.}
    \label{ml_mg_example}
\end{figure}

Figure \ref{ml_mg_example} shows a comparison of the results with the spline interpolation versus SR GAN interpolation and alternating spline-SR GAN interpolation ($N_{\rm{GAN}} = 1$). We plot the mean profiles and $1-\sigma$ bounds of the norm difference between the pressure grids at the current and previous iteration. The top row of the figure shows the result of using the SR operator at every iteration compared to the spline operator. The multigrid solver converges in 100 iterations (compared to $\sim 10^4$ iterations for a simple iterative solver). We find that the mean profile for the SR GAN operator performs just as well as the spline one. However, looking at $1-\sigma$ intervals of the ratio of the profiles reveals that there is a large variance in performance. Using the super resolution prolongation operator can result in faster convergence compared to bivariate spline operator, but for some fields, the spline interpolation can fare better. The bottom row shows the result of alternating the SR GAN and spline operators compared with using only the spline operator. There is consistently faster convergence using the alternating operators, as the ratio of SR and spline profiles is less than one. This improvement could be due to the fact that the SR GAN operator has been trained on real pressure field solutions and not random field, so it is more effective in conjunction with the spline operator than by itself. This result motivates a further exploration for multigrid parameters where SR GAN interpolation is particularly effective.

\subsection{Model parameters at fixed grid size}
\label{model_param_explore}
% N_GAN, Re, switching diff threshold value, N_smooth

Section \ref{motivatingresults} demonstrates the potential increased efficiency of the super resolution enhanced multigrid solver, although the efficiency may depend on choices in the model parameters. In order to explore the parameters, we fix the grid size to $N_{\rm{grid}} = 192$ and change $N_{\rm{smooth}}$, $N_{\rm{GAN}}$, $S_{\rm{thres}}$, and $\rm{Re}_{test}$.

%For these tests, we only solve 10 grids instead of the 100 grids in the previous section.

\begin{figure}

\begin{subfigure}{0.5\textwidth}
    \centering
    \includegraphics[width=\textwidth]{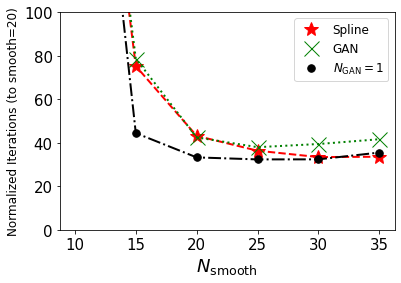}
    \caption{Plot of $N_{\rm{smooth}}$ vs. iterations for three interpolation types. The number of iterations are scaled by $N_{\rm{smooth}}$/20, accounting for the added/reduced computational effort in smoothing operations per iteration step above/below the fiducial value. All the data points for $N_{\rm{smooth}} = 10$ are outside the figure range, as the algorithm did not converge in the $N_{\rm{iter}} = 300$ max number of iterations we considered.}
    \label{nsmooth}
\end{subfigure}
\begin{subfigure}{0.5\textwidth}
        \centering
        \includegraphics[width=\textwidth]{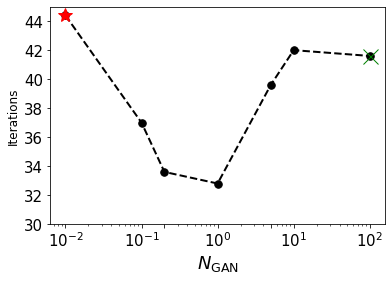}
        \caption{Plot of $N_{\rm{GAN}}$ vs. iterations. The endpoints of $<10^{-2}$ and $>10^{-2}$ represent the limiting cases of only spline and only SR GAN interpolation respectively. For these parameters, a value of $N_{\rm{GAN}} = 1$ performed the best.}
        \label{ngan}
\end{subfigure}
\vskip\baselineskip
\begin{subfigure}{0.5\textwidth}
    \centering
    \includegraphics[width=\textwidth]{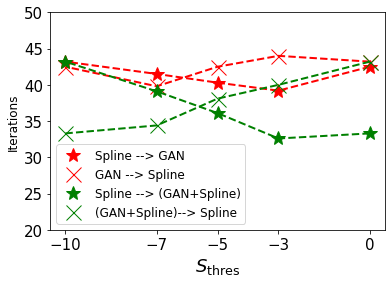}
    \caption{Plot of $S_{\rm{thres}}$ vs. iterations. Different permutations of before/after threshold interpolations are shown for spline, SR GAN, and $N_{\rm{GAN}} = 1$ interpolation types. The profiles denoted by `x' and `*' of a given color should be approximate reflections of each other.}
        \label{nthres}
\end{subfigure}
\begin{subfigure}{0.5\textwidth}
    \centering
    \includegraphics[width=\textwidth]{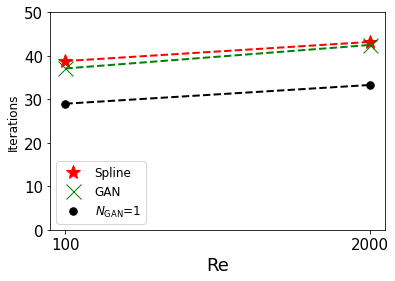}
    \caption{Plot of $\rm{Re}_{\rm{test}} = 200 and 2000$ vs. iterations for different interpolation cases.}
        \label{re_plot}
\end{subfigure}
\caption{Plots showing iterations to convergence for different parameters at a fixed $N_{\rm{grid}} = 192$.}

\end{figure}

Figure \ref{nsmooth} shows the number of iterations to convergence for different $N_{\rm{smooth}}$ values. We scale the number of iterations by a factor $N_{\rm{smooth}}$/20 in order to account for the differing total number of smoothing operations performed per multigrid iteration compared to the fiducial value of $N_{\rm{smooth}} = 20$. One smoothing operation is not quite equivalent to a multigrid iteration, but smoothing does help damp errors and thus get closer to convergence. The data points for $N_{\rm{smooth}} = 10$ are off the plot because no grid converged within $N_{\rm{iter}} = 300$ number of iterations. The plot also shows that for all cases after $N_{\rm{smooth}} = 20$ there is not an additional benefit to more smoothing in between iterations. The spline and SR GAN cases have similar average convergence for $N_{\rm{smooth}} < 25$, although the SR GAN profile is increasingly worse for $N_{\rm{smooth}} > 25$. The case of $N_{\rm{GAN}} = 1$ is the option with fastest convergence, except at $N_{\rm{smooth}} > 30$, where the spline case becomes comparable. It is also worth noting that the $N_{\rm{GAN}} = 1$ profile shows the least amount of variation in iterations to convergence, especially between  $N_{\rm{smooth}} = 15-20$.

Figure \ref{ngan} shows the number of iterations to convergence for different $N_{\rm{GAN}}$ values. The endpoints of the plot represent the limits of spline-only ($N_{\rm{GAN}} < 10^{-2}$) and SR GAN-only ($N_{\rm{GAN}} > 10^{2}$) interpolation (i.e. in theory, values should be 0 or $\infty$, but we have a finite number of iterations tested). The transition from spline-only to including some SR GAN interpolations causes a steep drop in iterations compared to going from SR GAN-only and including more spline interpolations. Interestingly, there was not much difference in cases with $N_{\rm{GAN}} > 10^{2}$ and $N_{\rm{GAN}} = 10$, as in, adding a spline-interpolation every 10 iterations did not improve results. The curve shows a minimum near $N_{\rm{GAN}}$ = 1, where  spline and SR GAN interpolations are alternated each iteration.

In addition to alternating spline and SR GAN interpolations, we can also set a cutoff based on how far the algorithm has progressed. This option allows us to apply the operators at particular parts of the calculation, so we use the operator that is most efficient for the current stage of solution \citep[][]{markidis2021old}. The parameter $S_{\rm{thres}}$ sets a transition switch in prolongation algorithm based on the norm of the difference in the pressure fields at pairs of iterations, $|p_i - p_{i+1}|$. Figure \ref{nthres} shows the number of iterations to convergence vs. $S_{\rm{thres}}$. Since convergence is defined as a tolerance $|p_i - p_{i+1}| = 10^{-10}$, $S_{\rm{thres}} = 10^{-10}$ represents the case where the there is no switch of the interpolation operator. The opposite case of $S_{\rm{thres}} = 10^{0}$ represents the case where there is an immediate switch from the first to the second interpolation choice. The plot shows that the iterations to convergence remains relatively similar when switching between spline and SR GAN operators. Following the `spline $\rightarrow$ GAN' profile, we see that switching to SR GAN earlier in the calculation (going left to right) reduces the number of iterations, with the best choice being $S_{\rm{thres}} = 10^{-3}$. The reverse `GAN $\rightarrow$ spline' profile shows a similar minimum iteration value at $S_{\rm{thres}} = 10^{-7}$. The similarity in the minimum value of iterations for algorithms starting with either spline or SR GAN could suggest that the important factor is not which operator starts the iterative process, but instead how often the SR GAN operator is used. The green profiles switching between spline and $N_{\rm{GAN}} = 1$ (`GAN+Spline') cases show that these options perform better than the previous cases, consistent with the finding in Figure \ref{ngan} that  $N_{\rm{GAN}} = 1$ is the most efficient choice for these particular parameters. The optimal choices of $S_{\rm{thres}}$ for these profiles are the ones that retain more of the $N_{\rm{GAN}} = 1$ operator choice.

We also examine the ability of the SR GAN model to represent data with different spatial scales than the model was trained on. The training data is generated with $\rm{Re}= 2000$. For test data, we also examined fluid fields with  $\rm{Re}= 100$. Figure \ref{re_plot} shows that there is not a large difference in convergence iterations and relative trends between interpolation choices with lower $\rm{Re}$. We do not test larger $\rm{Re}$ for these $N_{\rm{grid}}=192$ grids because the driving scales become smaller than the resolution scale. In Figure \ref{scaling_re_plot} we examine larger $\rm{Re}$ for $N_{\rm{grid}}=384$.

\subsection{Performance and scaling at different grid sizes}
\label{performance}

\begin{figure}
    \centering
    \includegraphics[width=0.75\textwidth]{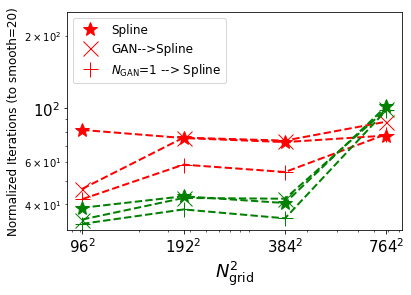}
    \caption{Profile of iteration to convergence vs. grid size for different parameter choices. For the examples with SR GAN interpolation, $S_{\rm{thres}} = 10^{-5}$. The data markers are of different choices of interpolation of spline-only, SR GAN $\rightarrow$ spline and $N_{\rm{GAN}} = 1$ $\rightarrow$ spline. The colors red and green are of choices of $N_{\rm{smooth}}$ = 15 and 20 respectively.}
    \label{scaling_nsmooth_plot}
\end{figure}

In Section \ref{model_param_explore} we fix the length of each dimension of the grid to $N_{\rm{grid}} = 192$ in order to explore the effects of various parameters. We now explore $N_{\rm{grid}}$ = 96, 192, 384, and 768 for a limited set of parameters. We examine both the cases of the Poisson source term being interpolated onto a different resolution grid from the fiducial case and the source being produced naturally by driving turbulence on a different-sized grid. The first case investigates the performance of the SR GAN operator at different samplings of a fixed spatial scale. The second case investigates the ability of the SR GAN operator to capture spatial structure on different scales than that of the training set. Due to the increased computational cost, we use 10 grids for the first case and 5 grids for the second case.

First, we assume the same underlying source term as before. The source term is interpolated from $N_{\rm{grid}} = 192$ to a different value via a two-dimensional spline interpolation. In these runs, the SR GAN multigrid algorithm converges to machine precision. The fiducial runs converged to a tolerance $<10^{-10}$, while the interpolated source term runs converged to tolerances between  $\sim 10^{-4}$ and $\sim 10^{-7}$. Even though the tolerance is larger, these grids do approach the true solution. It is likely that this behavior is caused by lack of physically realistic at spatial scales interpolated from the coarse source term to a finer grid with a spline method. The SR GAN operator will nonetheless add information at these newly added scales. Since the source term does not contain meaningful information at these scales, the additional information may not help the solution progress in convergence (see Figure \ref{SR_GAN_test} for a visual comparison of injected spatial scales by the different operators).

Figure \ref{scaling_nsmooth_plot} shows profiles with $S_{\rm{thres}} = 10^{-5}$, meaning that when the norm difference in grids between multigrid steps falls below $S_{\rm{thres}}$, the interpolation operator switches to only spline. We plot profiles for $N_{\rm{smooth}}$ = 15 and 20. We chose these two values because they spanned the region from no convergence to convergence in Figure \ref{nsmooth}. The number of iterations to convergence is multiplied by $(N_{\rm{smooth}}/20)$ in order to provide a fair comparison of the computational effort between choices. This normalization approximates one Gauss-Seidel smoothing iteration as one multigrid iteration, which is good enough for our rough comparison. We see that generally $N_{\rm{smooth}}=20$ (green) performs better than $N_{\rm{smooth}}=15$ (red), as Figure \ref{nsmooth} indicates. At $N_{\rm{grid}} = 192$ we see that the operator with $N_{\rm{GAN}}=1$ performs the best, as was found in Section \ref{model_param_explore}. The same trend is generally true for the other grid sizes. For the undersampled source term $N_{\rm{grid}} = 96$ grid, convergence is slightly better compared to $N_{\rm{grid}} = 192$. Interestingly, for $N_{\rm{smooth}}=15$ the SR GAN interpolators are significantly better than the spline. At the other grid size extreme of $N_{\rm{grid}} = 768$, convergence becomes worse for the SR GAN examples. In particular, the $N_{\rm{GAN}}= 1$ profile approaches the same convergence value as the spline one. The reduction in efficiency is even greater for $N_{\rm{smooth}}=20$, where all examples are comparable in iterations to $N_{\rm{smooth}}=15$.

\begin{figure}
\begin{subfigure}{0.5\textwidth}
    \centering
    \includegraphics[width=\textwidth]{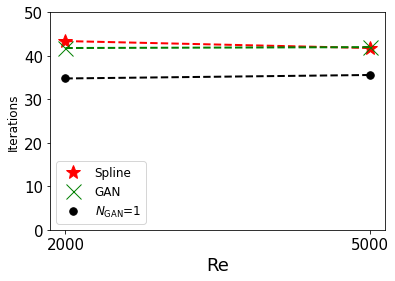}
    \caption{Plot similar to Figure \ref{re_plot}, except the source term grid is $N=384$. The colors are of choices of $\rm{Re}_{test}$ = 2000 and 5000. }
        \label{scaling_re_plot}
\end{subfigure}
\begin{subfigure} {0.5\textwidth}
    \centering
    \includegraphics[width=\textwidth]{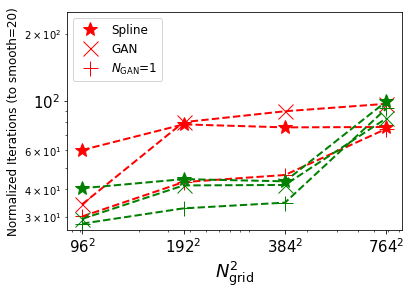}
    \caption{Plot similar to Figure \ref{scaling_nsmooth_plot}, except the parameter $S_{\rm{thres}}$ is not used. }
        \label{scaling_N_plot}
\end{subfigure}
    \caption{Convergence for different $\rm{Re}$ and grid size}
\end{figure}

Second, we solve the Poisson equation for larger grids using a source term derived from that particular grid size, instead of interpolation from $N_{\rm{grid}} = 192$. In Figure \ref{scaling_N_plot} we examine the ability of the SR GAN to capture higher frequency spatial scales than the data it was trained on. There is not a significant difference in convergence for three values of $N_{\rm{GAN}}$ at $\rm{Re} = 5000$. Figure \ref{scaling_N_plot} is similar to Figure \ref{scaling_nsmooth_plot}, except the source term is calculated from fluid velocity fields generated at the same grid resolution, instead of interpolated from the fiducial resolution of $N=192$. The general behavior of the profiles is similar to Figure \ref{scaling_nsmooth_plot}: convergence is better for $N_{\rm{smooth}} = 20$ compared to a value of 15, although performance suffers significantly at the highest resolution.

\section{Discussion}

% Interpretability
%%%% Discuss what the GAN is doing. Why is the GAN doing better? What are the small scales it's working better on?

It is important to note that the initial motivating results shown in Section \ref{motivatingresults} demonstrate that adding the SR GAN to the multigrid interpolation step produces an algorithm that converges reasonably well in terms of iterations to convergence. In other words, on average the algorithm still converges on the same grids that the spline operator does. This result demonstrates that the small-scale structure injected by the SR GAN interpolation operator, which is not present in the spline interpolation (as shown in Figure \ref{SR_GAN_test}), provides a realistic contribution to the pressure field solution as the multigrid algorithm progresses. Looking more closely at the large spread in performance of the SR GAN solver in Figure \ref{ml_mg_example}, we see that the SR GAN multigrid solver can perform better on some grids compared to the spline multigrid solver, although on other grids it performed worse.

Figure \ref{ml_mg_example} also shows that alternating use of the spline and SR GAN operators ($N_{\rm{GAN}} = 1$) produced a solver that consistently converges faster than the spline-based solver. Figure \ref{ngan} shows that $N_{\rm{GAN}} = 1$ seems to be the optimal choice. Also, the figure shows convergence is not symmetric around $N_{\rm{GAN}} = 1$; less SR GAN operations per spline operator ($N_{\rm{GAN}} < 1$) continued to perform relatively well as $N_{\rm{GAN}}$ became lower, while performance became worse more quickly as  $N_{\rm{GAN}}$ increased.

Along the same line of thought, we introduce the $S_{\rm{thres}}$ parameter, which provides a switch for changing the interpolation operator depending on how converged the solution is. We might expect that the SR GAN operator would perform better on grids that are closer to convergence, as these grids are closer to the GAN's learned representation of the training data. In contrast, we can also argue that the injection of smaller scale structure by the SR GAN at earlier times in the calculation could accelerate convergence. Figure \ref{nthres} shows that switching between spline and SR GAN operators (red lines) did not greatly improve convergence, but there is a preference towards including a mix of operators, as the minimum in convergence is found between the two extremes of $S_{\rm{thres}}$. Furthermore, the profiles are asymmetric in $S_{\rm{thres}}$, tending towards using the SR GAN operator for more of the calculation (lower $S_{\rm{thres}}$ starting with `GAN' and higher $S_{\rm{thres}}$ starting with `spline').The order of the choice of operators did not seem to matter much, so both of the opposing arguments earlier in the paragraph are unsupported. Similar asymmetry in the profiles occurs when switching between `spline' ($N_{\rm{GAN}} < 10^{-2}$ ) and `GAN + Spline' ($N_{\rm{GAN}} = 1$), except that $N_{\rm{GAN}} = 1$ is clearly preferred. One possible explanation for lack of optimal order of operator choices could be that each operator is more efficient at different spatial scales, as the power spectra in Figure \ref{SR_GAN_test_set} shows. Including both operators within the calculation targets more scales than only including either operator.

These finding follow the idea that ML applications to scientific domains could work best in tandem with existing methods, rather than as a whole replacement. Different operators can have benefits or drawbacks depending on the data at hand, so a mixed operator multigrid algorithm could have a more robust domain of effectiveness. \citet[][]{zhang2022hybrid} tested a different operator replacement scheme within multigrid than our work, and also found that a mixture of traditional and ML methods produced the best results. Our findings of additional hybrid algorithm parameters that did not end up providing significant benefit highlights the rich area of research in identifying and characterising parameters unique to the hybrid nature of the algorithm.

We also explored the transferability of the deep learning model. We tested grids with lower and higher $\rm{Re}$ than the training data in Figures \ref{re_plot} and \ref{scaling_re_plot}. Changing the data $\rm{Re}$ did not have a significant impact on the efficiency of the multigrid solver for any of the three cases tested. There was a slight improvement at $\rm{Re}=100$. The grid size limited our exploration of higher $\rm{Re}$. As mentioned in Section \ref{multigrid_section}, the multigrid method is applied recursively on grids of different residual scales because the residuals also satisfy the Poisson equation. Due to this self-similar property, we can expect that the SR GAN operator, which learned on grids with $\rm{Re}=2000$, will be effective on data of additional scales. This expectation is supported the that the performance of the SR GAN multigrid solver did not decrease significantly with changes in $\rm{Re}$. It would be interesting to test the solver on a system with a much larger range of spatial scales to see if this behavior remains.

Finally, we explore the behavior of the solvers with different grid sizes. We found a potential issue with the SR GAN multigrid solver in that it does not converge to machine tolerance when using a source term that has been oversampled (i.e. interpolated with a spline method to a finer grid), likely due to the lack of physically relevant frequencies in the added smaller scale regime. Nonetheless, it is still possible to achieve convergence to machine precision by switching via the $S_{\rm{thres}}$ parameter to the spline solver once the convergence stalls, as Figure \ref{scaling_nsmooth_plot} shows. In this case, the $N_{\rm{GAN}}=1$ solver still performs better than the traditional solver. If we instead generate source terms at each grid size, the SR GAN solver does not have issues converging to machine precision and the resulting profiles are quite similar. As we expect, the curves with $N_{\rm{smooth}}=20$ generally show better convergence compared to $N_{\rm{smooth}}=15$. Interestingly, with $N_{\rm{smooth}}=20$ performance suffers at $N=768$, but not at $N=384$. This is true even for the spline-based solver. Further investigation is needed to determine if even higher resolution grids have similar drops in performance.

\subsection{Caveats and further work}

\begin{itemize}

\item We only considered the number of multigrid iterations to convergence, and not the total run time. Our implementation of multigrid with only SR GAN is up to a factor of 2 slower in time. The scaling with time will depend on the underlying ML implementation (Tensorflow in our case). Generally the best performance is seen in $N_{\rm{GAN}} = 1$ (alternating use of spline and SR GAN) or near that value, so the time cost is up to a factor of 1.5 worse. More optimized implementations and advancements in ML codes could plausibly decrease this time cost further.

\item We chose a series of parameters to focus on. There are other parameters that we kept fixed into order to simply analysis. For example, we fix the parameters that split the domain up into kernels of size $n_{\rm{s}}^2$ with stride = 2. The overlapping sections of the kernels are ignored, whereas \citet[][]{margenberg2022} take the average value in overlapping sections.

\item Our implementation of multigrid is a two-level algorithm, in order to explore the effect of a SR GAN interpolation operator in a system that can struggles to converge. For the Poisson equation, a multilevel multigrid algorithm converges extremely quickly.

\item The Poisson equation provides a simple system, while still containing complex spatial scales, to evaluate our work. Solving a more challenging system, such as the Helmholtz equation \citep[e.g.][]{zhang2022hybrid, lerer2023multigrid}, is of interest.

\item Several successful GAN-based super resolution strategies have been applied to reconstruct high resolution fields from coarse-filtered fields in the literature \citep[e.g.][]{deng2019super, yousif2021high, kim2021unsupervised}. Our work can leverage any one of the state-of-the-art generative models and be applied in conjunction with the multigrid strategy, as an accelerator. Future work will also involve testing some of the newer generative models within the proposed framework. 

\end{itemize}

\section{Conclusions}

We implement a multigrid algorithm with two different interpolation operator choices (bivariate spline or deep learning super resolution GAN), and allow for several runtime choices of the two options. We focus on solving the two-dimensional Poisson equation from the pressure-Poisson formulation of the Navier-Stokes equations.  We explore a set of key parameters for this algorithm, including the number of smoothing operations in between iterations, the number of GAN operators per spline operator, a threshold switch between operators depending on how converged the solution is, and the Reynolds number of test data. In our investigation, we conclude the following:
\begin{enumerate}
\item We demonstrate that replacing the interpolation step of the multigrid solver with a SR GAN model on average performs similarly to a solver with traditional spline interpolation. The variance in performance with individual grids motivates exploration of additional multigrid parameters. 

\item In our parameter study, we find that generally algorithms with mixed operator choice (applying different operators at particular intervals or closeness to solution convergence) converge with less iterations compared to algorithms with purely traditional spline or SR GAN methods. Performance varies with particular parameters unique to the hybrid algorithm. More detailed exploration of this category of parameters in any hybrid ML/tradition numerical scheme should be performed.

\item  We find that the SR GAN-based multigrid solver performs well on grids that are lower/higher resolution by a factor of two. However, there is a noticeable degradation in performance at the highest resolution grid tested, a factor of four increase from the fiducial resolution. 

\end{enumerate}

\section{Declaration of competing interest}
The authors declare that they have no known competing financial interests or personal relationships that could have appeared to influence the work report in this paper.

\section{Data availability}
Data will be made available on reasonable request.

\section{Acknowledgements}
We thank Natalie Klein for useful conversations. F.H. acknowledges support from the 2020 XCP Computational Physics workshop at Los Alamos National Laboratory, the NASA FINESST fellowship (grant number 80NSSC20K1541), and the University of Michigan Rackham Predoctoral fellowship.

Computational resources for this work were provided by Los Alamos National Laboratory.

The contributions from GDP were performed under the auspices of the U.S. Department of Energy by Lawrence Livermore National Laboratory under Contract DE-AC52-07NA27344.

%%%%%%%%%%%%%%%%%%%%%%%%%%
%\bibliographystyle{elsarticle}
\bibliographystyle{abbrvnat}
\bibliography{main}

\begin{thebibliography}{33}
\providecommand{\natexlab}[1]{#1}
\providecommand{\url}[1]{\texttt{#1}}
\expandafter\ifx\csname urlstyle\endcsname\relax
  \providecommand{\doi}[1]{doi: #1}\else
  \providecommand{\doi}{doi: \begingroup \urlstyle{rm}\Url}\fi

\bibitem[Al-Saffar et~al.(2017)Al-Saffar, Tao, and Talab]{al2017review}
A.~A.~M. Al-Saffar, H.~Tao, and M.~A. Talab.
\newblock Review of deep convolution neural network in image classification.
\newblock In \emph{2017 International Conference on Radar, Antenna, Microwave,
  Electronics, and Telecommunications (ICRAMET)}, pages 26--31. IEEE, 2017.

\bibitem[Birla()]{birla2018}
D.~Birla.
\newblock URL \url{https://github.com/deepak112/Keras-SRGAN}.

\bibitem[Creswell et~al.(2018)Creswell, White, Dumoulin, Arulkumaran, Sengupta,
  and Bharath]{creswell2018generative}
A.~Creswell, T.~White, V.~Dumoulin, K.~Arulkumaran, B.~Sengupta, and A.~A.
  Bharath.
\newblock Generative adversarial networks: An overview.
\newblock \emph{IEEE Signal Processing Magazine}, 35\penalty0 (1):\penalty0
  53--65, 2018.

\bibitem[Deng et~al.(2019)Deng, He, Liu, and Kim]{deng2019super}
Z.~Deng, C.~He, Y.~Liu, and K.~C. Kim.
\newblock Super-resolution reconstruction of turbulent velocity fields using a
  generative adversarial network-based artificial intelligence framework.
\newblock \emph{Physics of Fluids}, 31\penalty0 (12), 2019.

\bibitem[Farsiu et~al.(2004)Farsiu, Robinson, Elad, and
  Milanfar]{farsiu2004advances}
S.~Farsiu, D.~Robinson, M.~Elad, and P.~Milanfar.
\newblock Advances and challenges in super-resolution.
\newblock \emph{International Journal of Imaging Systems and Technology},
  14\penalty0 (2):\penalty0 47--57, 2004.

\bibitem[Goodfellow et~al.(2014)Goodfellow, Pouget-Abadie, Mirza, Xu,
  Ward-Farley, Ozair, Courville, and Bengio]{goodfellow2014}
I.~J. Goodfellow, J.~Pouget-Abadie, M.~Mirza, B.~Xu, D.~Ward-Farley, S.~Ozair,
  A.~Courville, and Y.~Bengio.
\newblock Generative adverserial nets.
\newblock \emph{Adv. Neural Procs. Sys.}, 27:\penalty0 2672 -- 2680, 2014.

\bibitem[Harshvardhan et~al.(2020)Harshvardhan, Gourisaria, Pandey, and
  Rautaray]{harshvardhan2020comprehensive}
G.~Harshvardhan, M.~K. Gourisaria, M.~Pandey, and S.~S. Rautaray.
\newblock A comprehensive survey and analysis of generative models in machine
  learning.
\newblock \emph{Computer Science Review}, 38:\penalty0 100285, 2020.

\bibitem[He et~al.(2019)He, Li, Feng, Ho, Ravanbakhsh, Chen, and
  P{\'o}czos]{he2019}
S.~He, Y.~Li, Y.~Feng, S.~Ho, S.~Ravanbakhsh, W.~Chen, and B.~P{\'o}czos.
\newblock Learning to predict the cosmological structure formation.
\newblock \emph{Proceedings of the National Academy of Sciences}, 116\penalty0
  (28):\penalty0 13825--13832, 2019.

\bibitem[Holguin et~al.(2021)Holguin, Sidharth, and
  Portwood]{holguin2021multigrid}
F.~Holguin, G.~Sidharth, and G.~Portwood.
\newblock Multigrid solver with super-resolved interpolation.
\newblock \emph{arXiv preprint arXiv:2105.01739}, 2021.

\bibitem[Jabbar et~al.(2021)Jabbar, Li, and Omar]{jabbar2021survey}
A.~Jabbar, X.~Li, and B.~Omar.
\newblock A survey on generative adversarial networks: Variants, applications,
  and training.
\newblock \emph{ACM Computing Surveys (CSUR)}, 54\penalty0 (8):\penalty0 1--49,
  2021.

\bibitem[Kim et~al.(2021)Kim, Kim, Won, and Lee]{kim2021unsupervised}
H.~Kim, J.~Kim, S.~Won, and C.~Lee.
\newblock Unsupervised deep learning for super-resolution reconstruction of
  turbulence.
\newblock \emph{Journal of Fluid Mechanics}, 910:\penalty0 A29, 2021.

\bibitem[Kobyzev et~al.(2020)Kobyzev, Prince, and
  Brubaker]{kobyzev2020normalizing}
I.~Kobyzev, S.~Prince, and M.~Brubaker.
\newblock Normalizing flows: An introduction and review of current methods.
\newblock \emph{IEEE Transactions on Pattern Analysis and Machine
  Intelligence}, 2020.

\bibitem[LeCun and Hinton(2015)]{lecunhinton2015}
Y.~LeCun and G.~Hinton.
\newblock Deep learning.
\newblock \emph{Nature}, 521:\penalty0 436 -- 444, 2015.

\bibitem[Ledig et~al.(2017)Ledig, Theis, Husz{\'a}r, Caballero, Cunningham,
  Acosta, Aitken, Tejani, Totz, Wang, et~al.]{ledig2017}
C.~Ledig, L.~Theis, F.~Husz{\'a}r, J.~Caballero, A.~Cunningham, A.~Acosta,
  A.~Aitken, A.~Tejani, J.~Totz, Z.~Wang, et~al.
\newblock Photo-realistic single image super-resolution using a generative
  adversarial network.
\newblock pages 4681--4690, 2017.

\bibitem[Lerer et~al.(2023)Lerer, Ben-Yair, and Treister]{lerer2023multigrid}
B.~Lerer, I.~Ben-Yair, and E.~Treister.
\newblock Multigrid-augmented deep learning for the helmholtz equation: Better
  scalability with compact implicit layers.
\newblock \emph{arXiv preprint arXiv:2306.17486}, 2023.

\bibitem[Margenberg et~al.(2022)Margenberg, Hartmann, Lessig, and
  Richter]{margenberg2022}
N.~Margenberg, D.~Hartmann, C.~Lessig, and T.~Richter.
\newblock A neural network multigrid solver for the navier-stokes equations.
\newblock \emph{Journal of Computational Physics}, page 110983, 2022.

\bibitem[Markidis(2021)]{markidis2021old}
S.~Markidis.
\newblock The old and the new: Can physics-informed deep-learning replace
  traditional linear solvers?
\newblock \emph{Frontiers in big Data}, page~92, 2021.

\bibitem[Minaee et~al.(2021)Minaee, Kalchbrenner, Cambria, Nikzad, Chenaghlu,
  and Gao]{minaee2021deep}
S.~Minaee, N.~Kalchbrenner, E.~Cambria, N.~Nikzad, M.~Chenaghlu, and J.~Gao.
\newblock Deep learning--based text classification: A comprehensive review.
\newblock \emph{ACM Computing Surveys (CSUR)}, 54\penalty0 (3):\penalty0 1--40,
  2021.

\bibitem[Park et~al.(2003)Park, Park, and Kang]{park2003super}
S.~C. Park, M.~K. Park, and M.~G. Kang.
\newblock Super-resolution image reconstruction: a technical overview.
\newblock \emph{IEEE signal processing magazine}, 20\penalty0 (3):\penalty0
  21--36, 2003.

\bibitem[Peurifoy et~al.(2018)Peurifoy, Shen, Jing, Yang, Cano-Renteria,
  DeLacy, Joannopoulos, Tegmark, and Solja{\v{c}}i{\'c}]{peurifoy2018}
J.~Peurifoy, Y.~Shen, L.~Jing, Y.~Yang, F.~Cano-Renteria, B.~G. DeLacy, J.~D.
  Joannopoulos, M.~Tegmark, and M.~Solja{\v{c}}i{\'c}.
\newblock Nanophotonic particle simulation and inverse design using artificial
  neural networks.
\newblock \emph{Science advances}, 4\penalty0 (6):\penalty0 eaar4206, 2018.

\bibitem[Portwood et~al.(2021)Portwood, Nadiga, Saenz, and Livescu]{portwood21}
G.~D. Portwood, B.~T. Nadiga, J.~A. Saenz, and D.~Livescu.
\newblock Interpreting neural network models of residual scalar flux.
\newblock \emph{Journal of Fluid Mechanics}, 907:\penalty0 A23, 2021.

\bibitem[Quartapelle(2013)]{quartapelle2013numerical}
L.~Quartapelle.
\newblock \emph{Numerical solution of the incompressible Navier-Stokes
  equations}, volume 113.
\newblock Birkh{\"a}user, 2013.

\bibitem[Raissi et~al.(2019)Raissi, Perdikaris, and
  Karniadakis]{raissi2019physics}
M.~Raissi, P.~Perdikaris, and G.~E. Karniadakis.
\newblock Physics-informed neural networks: A deep learning framework for
  solving forward and inverse problems involving nonlinear partial differential
  equations.
\newblock \emph{Journal of Computational Physics}, 378:\penalty0 686--707,
  2019.

\bibitem[Ranade et~al.(2020)Ranade, Hill, and Pathak]{ranade2020}
R.~Ranade, C.~Hill, and J.~Pathak.
\newblock Discretizationnet: A machine-learning based solver for navier-stokes
  equations using finite volume discretization.
\newblock \emph{arXiv preprint arXiv:2005.08357}, 2020.

\bibitem[Sanchez-Gonzalez et~al.(2020)Sanchez-Gonzalez, Godwin, Pfaff, Ying,
  Leskovec, and Battaglia]{sanchez2020}
A.~Sanchez-Gonzalez, J.~Godwin, T.~Pfaff, R.~Ying, J.~Leskovec, and
  P.~Battaglia.
\newblock Learning to simulate complex physics with graph networks.
\newblock In \emph{International Conference on Machine Learning}, pages
  8459--8468. PMLR, 2020.

\bibitem[Sidharth(2020)]{gs2020multiscale}
G.~S. Sidharth.
\newblock A multiscale subgrid decomposition.
\newblock \emph{AIAA Scitech 2020 Forum}, page 0820, 2020.

\bibitem[Trottenberg et~al.(2001)Trottenberg, Oosterlee, and
  Schuller]{trottenberg2001}
U.~Trottenberg, C.~W. Oosterlee, and A.~Schuller.
\newblock \emph{Multigrid}.
\newblock 2001.

\bibitem[Virtanen et~al.(2020)Virtanen, Gommers, Oliphant, Haberland, Reddy,
  Cournapeau, Burovski, Peterson, Weckesser, Bright, et~al.]{virtanen2020}
P.~Virtanen, R.~Gommers, T.~E. Oliphant, M.~Haberland, T.~Reddy, D.~Cournapeau,
  E.~Burovski, P.~Peterson, W.~Weckesser, J.~Bright, et~al.
\newblock Scipy 1.0: fundamental algorithms for scientific computing in python.
\newblock \emph{Nature methods}, 17\penalty0 (3):\penalty0 261--272, 2020.

\bibitem[Wang et~al.(2020)Wang, Kashinath, Mustafa, Albert, and Yu]{wang20}
R.~Wang, K.~Kashinath, M.~Mustafa, A.~Albert, and R.~Yu.
\newblock Towards physics-informed deep learning for turbulent flow prediction.
\newblock In \emph{Proceedings of the 26th ACM SIGKDD International Conference
  on Knowledge Discovery \& Data Mining}, pages 1457--1466, 2020.

\bibitem[Wei et~al.(2018)Wei, Jiang, and Chen]{wei2018}
Q.~Wei, Y.~Jiang, and J.~Z. Chen.
\newblock Machine-learning solver for modified diffusion equations.
\newblock \emph{Physical Review E}, 98\penalty0 (5):\penalty0 053304, 2018.

\bibitem[Yang et~al.(2019)Yang, Zhang, Tian, Wang, Xue, and Liao]{yang2019deep}
W.~Yang, X.~Zhang, Y.~Tian, W.~Wang, J.-H. Xue, and Q.~Liao.
\newblock Deep learning for single image super-resolution: A brief review.
\newblock \emph{IEEE Transactions on Multimedia}, 21\penalty0 (12):\penalty0
  3106--3121, 2019.

\bibitem[Yousif et~al.(2021)Yousif, Yu, and Lim]{yousif2021high}
M.~Z. Yousif, L.~Yu, and H.-C. Lim.
\newblock High-fidelity reconstruction of turbulent flow from spatially limited
  data using enhanced super-resolution generative adversarial network.
\newblock \emph{Physics of Fluids}, 33\penalty0 (12), 2021.

\bibitem[Zhang et~al.(2022)Zhang, Kahana, Turkel, Ranade, Pathak, and
  Karniadakis]{zhang2022hybrid}
E.~Zhang, A.~Kahana, E.~Turkel, R.~Ranade, J.~Pathak, and G.~E. Karniadakis.
\newblock A hybrid iterative numerical transferable solver (hints) for pdes
  based on deep operator network and relaxation methods.
\newblock \emph{arXiv preprint arXiv:2208.13273}, 2022.

\end{thebibliography}
%%%%%%%%%%%%%%%%%%%%%%%%%%

\appendix

\section{Algorithms}
\label{alg_appendix}

% see https://www.overleaf.com/learn/latex/Algorithms

\begin{algorithm}[H]
    \caption{\emph{Linear  normalization}: This function maps dimensional data $\textbf{p}$ to normalized space. }
    \textbf{Function} Linear-Normalize($\mathbf{p}$, $u$, $v$, $a$, $b$)
    \begin{algorithmic}
        \State $\textbf{q} \gets (\mathbf{p} - u) \left( \frac{ b-a}{v-u} \right)$ 
        \State \textbf{return} $\textbf{q}$
    \end{algorithmic}
\end{algorithm}

\begin{algorithm}[H]
    \caption{\emph{Linear  denormalization}: This function maps back non-dimensional data $\textbf{p}$ back to its dimensional values. }
    \textbf{Function} Linear-Deormalize($\mathbf{q}$, $u$, $v$, $a$, $b$)
    \begin{algorithmic}
        \State $\textbf{p} \gets (\mathbf{q} - a) \left( \frac{v-u}{b-a} \right) + u$
        \State \textbf{return} $\textbf{p}$
    \end{algorithmic}
\end{algorithm}

\begin{algorithm}[H]
    \caption{\emph{Symmetric  normalization}: This function maps dimensional data using a logarithmic normalization  to preserve the sign of the data}
    \textbf{Function} Normalize($\mathbf{p}$, $p_{\textrm{max}}$, $p_{\textrm{min}}$ )
    \begin{algorithmic}
        \State $\mathbf{p}', p_{\textrm{min}}',p_{\textrm{max}}' \gets \textrm{log10}(|\mathbf{p}|), \textrm{log10}(p_{\textrm{min}}), \textrm{log10}(p_{\textrm{max}})$
        \State $\mathbf{q} \gets \textrm{Linear-Normalize}(\mathbf{p}', p_{\textrm{min}}' , p_{\textrm{max}}' ,0,1  ) \circ \rm{sgn}(\mathbf{p}) $
        \State \textbf{return} $\mathbf{q}, \rm{sgn}(\mathbf{p})$
    \end{algorithmic}
\end{algorithm}

\begin{algorithm}[H]
    \caption{\emph{Symmetric  denormalization}: This function maps back logarithmicaly normalized data to dimensional data.}
    \textbf{Function} Denormalize($\mathbf{q}$, $\mathbf{p}_{\textrm{sign}}$ $p_{\textrm{max}}$, $p_{\textrm{min}}$)
    \begin{algorithmic}
        \State $p_{\textrm{min}}',p_{\textrm{max}}' \gets \textrm{log10}(p_{\textrm{min}}), \textrm{log10}(p_{\textrm{max}}) $
        \State $\mathbf{q} \gets \textrm{Linear-Denormalize}(\mathbf{q}, p_{\textrm{min}}', p_{\textrm{max}}' ,0,1  )   $
        \State $\mathbf{return} \ \textrm{power}(10,\ \mathbf{q}) \circ \mathbf{p}_{\textrm{sign}}  $
    \end{algorithmic}
\end{algorithm}

\begin{algorithm}[H]
    \caption{\emph{Super resolution Prolongation}: This function describes the use of the GAN architecture in the prolongation operator of the Multigrid algorithm \ref{multigrid_algorithm}}
    \textbf{Function} SR-Prolongation($\mathbf{p}^h$)
    \begin{algorithmic}
        \State $\mathbf{q}^H \gets \mathbf{0}^H$
        \State $ \mathbf{q}^{h}, \mathbf{p}^h_{\rm{sign}} \gets \textrm{Normalize}(\mathbf{p}, p_{\textrm{max}}, p_{\textrm{min}} )$
        \State $(i_{\textrm{w}}, j_{\textrm{w}})  \gets 0$
        \While{$i_{\textrm{w}} < (N_{\textrm{grid}}-n_s) $}
            \While{$j_{\textrm{w}} < (N_{\textrm{grid}}-n_s) $}
                \State $\textrm{window indices} \gets (i_{\textrm{w}}:i_{\textrm{w}} + n_s, j_{\textrm{w}}:j_{\textrm{w}}+n_s )$

                \State $\textrm{central window indices} \gets (i_{\textrm{w}}+2:i_{\textrm{w}} + 4, j_{\textrm{w}}+2:j_{\textrm{w}}+4 )$
                
                \State $\mathbf{q}^{h}_{\textrm{w}} \gets  \mathbf{q}^h[\textrm{window indices}]$ 
                \State $H \gets 2 n_{\textrm{up}} h $
                \State $\mathbf{q}_{\textrm{w}}^i \gets \mathbf{q}_{\textrm{w}}^h$
                \For{2 iterations}
                \State $\mathbf{q}^{i}_{\textrm{w}} \gets $ Super-Resolution GAN($\mathbf{q}_{\textrm{w}}^i$) 
                \EndFor
                \State $\mathbf{q}_{\textrm{w}}^H \gets \mathbf{q}_{\textrm{w}}^i$
                \State $\textrm{up-sampled central window indices} \gets 2n_{\textrm{up}} \times \textrm{central window indices}$  
                \State $\mathbf{q}^H[\textrm{up-sampled central window indices} ] \gets \mathbf{q}^{H}_{\textrm{w}}$
                \State Fill any non-central, edge elements in $\mathbf{p}^H$ with data from  ${q}^{H}_{\textrm{w}}$
                \State $j_{\textrm{w}} \gets j_{\textrm{w}} + \textrm{stride}$
            \EndWhile
            \State $i_{\textrm{w}} \gets i_{\textrm{w}} + \textrm{stride}$
        \EndWhile
    \State $\mathbf{p}^H \gets \textrm{Denormalize}(\mathbf{q}^{H} , \mathbf{p}^h_{\rm{sign}},p_{\textrm{max}}, p_{\textrm{min}} )$
    \end{algorithmic}
\end{algorithm}

% V-cycle

\begin{algorithm}[H]
    \caption{Multigrid algorithm with V-Cycle structure}
    \textbf{Function} V-CYCLE($\mathbf{p, r}$)
    \begin{algorithmic}
        \State Relax $N_{\rm{iter}}$ times via Gauss-Seidel method
        \State $\mathbf{\delta}^h \gets \mathbf{0}$
        \If{h is coarsest level}
            \State Relax 10 times via Gauss-Seidel method to solve $\mathbf{\nabla \delta^{h} = -r}$
            \State \textbf{return} $\mathbf{\delta^{h}}$
        \EndIf
        \State $\mathbf{r} \gets \mathbf{f - A p}$
        \State $\mathbf{r}^{2h}, \mathbf{\delta}^{2h} \gets \textrm{Restriction}( \mathbf{r}, \delta)$
        \State $\mathbf{\delta}^{2h}  \gets \textrm{V-cycle}(\mathbf{r}, \mathbf{\delta }$)
        \State Solve $\mathbf{\nabla \delta^{h} = -r}$

        \If{$N_{\rm{GAN}} < 1$}
            \State $N = 1/N_{\rm{GAN}}$ 
            \State First Prolongation $\gets$ Bivariate Spline Prolongation
            \State Second Prolongation $\gets$ SR-Prolongation
            \State switch
        \Else
            \State $N = N_{\rm{GAN}}$ 
            \State First Prolongation $\gets$ SR-Prolongation
            \State Second Prolongation $\gets$ Bivariate Spline Prolongation
        \EndIf
        % ($\textrm{log}_{10}||\mathbf{r}||\leq {S_{\textrm{thres}} } $ ) }

        \State Switch $\gets$ false
        \If {norm of difference between previous iteration$\mathbf{p}_{i-1}$ and current $\mathbf{p}_{i} \leq S_{\textrm{thres}}$ }
            \State Switch $\gets$ true
        \EndIf
        
        \If { (current iteration mod $N$ is not equal to 0) and (not switch)  }
            \State $\mathbf{\delta}^h \gets \textrm{First Prolongation}\left(\mathbf{\delta}^{2h} \right)$
        \Else
                \quad\State $\mathbf{\delta}^h \gets \textrm{Second Prolongation}(\mathbf{\delta}^{2h} )$
        \EndIf
        \State \textbf{return} $ \mathbf{p + \delta}^h$
    \end{algorithmic}
    \label{multigrid_algorithm}
\end{algorithm}

\end{document}